\documentclass[11pt]{article}

\usepackage{amsmath,amsfonts,amssymb} 
\usepackage[english]{babel}                        
\usepackage[active]{srcltx}                         
\usepackage{geometry}
\usepackage[T1]{fontenc}			                                

\usepackage[english]{babel}                        
\usepackage[T1]{fontenc}			   

\usepackage{theorem}{\theorembodyfont{\slshape}
  \newtheorem{theorem}{Theorem}
  \newtheorem{lemma}{Lemma}
}{\theorembodyfont{\upshape}
  \newtheorem{remark}{Remark} 
}
\def\Proof{{\smallskip\noindent{\em Proof. }}}     
\def\endProof{{\hfill$\Box$\medskip\noindent}}     

\def\endProof{{\hfill$\Box$}}
\newcommand\R{{\mathbb{R}}}
\renewcommand\P{{\mathbb{P}}}
\newcommand\N{{\mathbb{N}}}

\renewcommand\S{{\mathbb{S}}}

\title{Concentration-diffusion effects in viscous incompressible flows\footnote{
{\it Keywords:\/} Navier-Stokes, asymptotic profiles, asymptotic behavior, far-field, 
flow map, analyticity,
symmetry, spatial spreading, lower bound, pointwise estimates,
decay rate, singularities.\hfill\break
{\it 2000 Mathematics subject classification:\/} 76D05, 35Q30
}}
\author{
Lorenzo Brandolese\footnote
{
Universit\'e de Lyon,
Universit\'e Lyon 1,
CNRS UMR 5208 Institut Camille Jordan,
43, bd. du 11 novembre 1918,
F - 69622 Villeurbanne Cedex,
France. E-mail:
{\tt brandolese@math.univ-lyon1.fr}
}
}

\begin{document}

\maketitle

\begin{abstract}
Given a finite sequence of times $0<t_1<\dots<t_N$\,,
we construct an example of a smooth solution of the free nonstationnary Navier--Stokes equations in $\R^d$,
$d=2,3$, 
such that:
(i) The velocity field $u(x,t)$ is spatially poorly localized
at the beginning of the evolution
but tends to concentrate until, as  the time~$t$ approaches $t_1$, it becomes well-localized.
(ii)
Then $u$ spreads out again after $t_1$, and
such concentration-diffusion phenomena are later reproduced near the instants
$t_2$, $t_3$,~\dots
\end{abstract}

\section{Introduction}

One of the most important questions in mathematical Fluid Mechanics,
which is still far from being understood, is to know whether
a finite energy, and initially smooth, nonstationnary Navier--Stokes flow
will always remain regular during its evolution,
or can become turbulent in finite time.

As a first step toward the understanding of possible blow-up mechanisms,
it is interesting to exhibit examples of smooth and decaying initial data
such that, even if the corresponding solutions
remain regular for all time, ``something strange'' happens
around a given point $(x_0,t_0)$ in space-time.
This is the goal of the present paper.

Our main result
is the construction of a class of (smooth) solutions
to the incompressible
Navier--Stokes equations
such that,
in the absence of any external forces,
the motion of the fluid particles tends to be more concentrated
around $x_0$, as the time~$t$
approaches~$t_0$.
This corresponds precisely to the qualitative behavior that one would
expect in the presence of a singularity,
even though  such ``concentration of the motion'' 
is not strong enough to imply their formation.

\subsection{Statement of the main result}

For a fluid filling the whole space~$\R^d$, $d\ge2$, the Navier--Stokes system
can be written as
\begin{equation*}
 \left\{
\begin{aligned}
 &\partial_t u+\P\nabla\cdot (u\otimes u)=\Delta u \\
 &\nabla \cdot u=0\\
&u(x,0)=a(x),
\end{aligned}
\right.
\end{equation*}
where $u=(u_1,\ldots,u_d)$, $a$ is a divergence-free vector field in $\R^d$
and $\P={\rm Id}-\nabla\Delta^{-1}\hbox{div}$ is the Leray-Hopf projector.

Because of the presence of the non-local operator~$\P$ 
a velocity field that is spatially  well localized (say, rapidly decaying as $|x|\to\infty$)
at the beginning of the evolution, in general, will immediately spread out.
A sharp description of this phenomenon
is provided by the two estimates~\eqref{ulbe} below.
In order to rule out the case of somewhat pathological flows (such as
two dimensional flows with radial vorticity, or the three-dimensional flows described in~\cite{Bra04iii},
which behave quite differently
as $|x|\to\infty$ if compared with generic solutions),
we will restrict our attention to data satisfying the following mild~\emph{non-symmetry} assumption
(for $j,k=1,\ldots,d$):
\begin{equation}
\label{nonor}
\exists\,j\not=k\colon\; \int (a_ja_k)(x)\,dx\not=0 ,\qquad
\hbox{or}\quad\int a_j(x)^2\,dx\not=\int a_k(x)^2\,dx.
\end{equation}
Then, for sufficiently fast decaying data, we have, for $|x|\ge \frac{C}{\sqrt t}$   (see~\cite{BraV07}):
\begin{equation}
 \label{ulbe}
\eta_1(t)|x|^{-d-1}\le |u(x,t)|\le \eta_2(t)|x|^{-d-1}, \qquad 
\hbox{$  \frac{x}{|x|}\in \S^{d-1}\backslash\Sigma_a$},
\end{equation}
these estimates being valid during  a small
 time interval $t\in(0,t_1)$.
Here $C,t_1>0$
and $\eta_1$  and $\eta_2$ are positive functions, independent on $x$, 
behaving like $\sim c_j\,t$ as $t\to0$ ($j=1,2$).
Moreover, $\S^{d-1}$ denotes the unit sphere in $\R^d$ and
the  subset $\Sigma_a$ of $\S^{d-1}$
represents the directions along which the lower bound may fail to hold:
the result of~\cite{BraV07} tells us that
$\Sigma_a$ can be taken of arbitrarily small surface measure on the sphere.
In other words, \emph{the lower bound holds  true in quasi-all directions\/},
whereas the upper bound is valid along all directions.

Moreover, the upper bound will hold  during the whole lifetime of the
strong solution~$u$ (see~\cite{Miy02}, \cite{Vig05}), whereas the lower bound
is valid, {\it a priori\/}, only during a very short time interval.
The main reason for this is that the
matrix $(\int u_ju_k(x,t)\,dx)$ is non-invariant during
the Navier--Stokes evolution, in a such way that even if the datum
satisfies~\eqref{nonor}, it cannot be excluded that at later times
the solution features some kind of creation of symmetry,
yielding to a better spatial localization and, after $t>t_1$,
to an improved decay as $|x|\to\infty$.
(We refer {\it e.g.\/} to \cite{BraM02}, \cite{Bra04iii},  \cite{HeM03}, \cite{HeMiy06}, for
the connection between the symmetry and the decay of solutions).

\medskip
The purpose of this paper is to show that this indeed can happen.
We construct an example of a solution of the Navier--Stokes
equations, with datum $a\in\mathcal{S}(\R^d)$ (the Schwartz  class), $d=2,3$,
such that the lower bound 
\begin{equation}
\label{lb2}
 |u(x,t)|\ge \eta_1(t)|x|^{-d-1},
\end{equation}
 holds in some interval $(0,t_1)$, but then brakes down at $t_1$, where
a \emph{stronger upper bound\/}  can be established.
This means that the motion of the fluid concentrates around the origin
at such instant.
Then the lower bound~\eqref{lb2} will hold true again after $t_1$,
until it will break down once more at a time $t_2>t_1$.
This diffusion-concentration effect can be repeated an arbitrarily
large number of times. 

More precisely, we will prove the following theorem.


\begin{theorem}
\label{theorem1}
Let $d=2,3$,  let $0<t_1<\cdots<t_N$ be a finite sequence and $\epsilon>0$.
Then there exist
a divergence-free vector field $a\in\mathcal{S}(\R^d)$ and
two sequences $(t_1',\ldots,t_N')$ and $(t_1^*,\ldots,t_N^*)$
such that the corresponding unique strong
solution $u(x,t)$ of the Navier--Stokes system
satisfies, for all $i=1,\ldots,N$  and all $|x|$ large enough,
the pointwise the lower bound
\begin{equation*}
|u(x, t'_i)|\ge c_\omega|x|^{-d-1}, 
\end{equation*}
and the stronger upper bound
\begin{equation*}
|u(x,t^*_i)|\le C|x|^{-d-2},
\end{equation*}
for a constant $C>0$ independent on~$x$
and a constant $c_\omega$ independent on $|x|$,
but possibly dependent on the projection~$\omega=\frac{x}{|x|}$
of $x$ on the sphere, and such that $c_\omega>0$ for a.e. $\omega\in\S^{d-1}$.
Moreover,
$t'_i$ and $t^*_i$ can be taken arbitrarily close to $t_i$:
\begin{equation*}
 |t'_i-t_i|<\epsilon \qquad\hbox{and}\qquad |t^*_i-t_i|<\epsilon,
\qquad\hbox{for $i=1,\ldots,N$}.
\end{equation*}
\end{theorem}

\begin{remark}
The initial datum can be chosen of the form
$a=\hbox{curl}(\psi)$,
where $\psi$ is a linear combination of dilated and modulated of a single function
(or vector field, if $d=3$)~$\phi\in\mathcal{S}(\R^d)$,
with compactly supported Fourier transform.
\end{remark}

Roughly speaking, our construction works as follows:
we look for an initial datum of the form 
$a=\hbox{curl}(\psi)$, where
\begin{equation*}
\label{psifo}
\psi(x)=\sum_{j=1}^{d(N+1)} \lambda_j\,\delta^{d/2}\phi(\delta x)\cos(\alpha_j\cdot x).
\end{equation*}
The unknown vector~$\boldsymbol{\alpha}=(\alpha_1,\ldots,\alpha_{d(N+1)})\in\R^{d^2(N+1)}$
of all the phases $\alpha_j\in\R^d$
will be assumed to belong to a suitable subspace $V\subset \R^{d^2(N+1)}$ of dimension~$d(N+1)$
in order to ensure, {\it a priori\/}, some nice geometrical properties of the flow.
Such geometric properties consist of a kind of rotational symmetry,
similar to that considered in~\cite{BraM02}, but less stringent.
In this way, the problem can be reduced to the study of the zeros
of the real function
$$t\mapsto\int_0^t\!\!\int (u_1u_2)(x,s)\,dx\,ds.$$
By an analyticity argument, this  in turn is reduced to
the study of the sign of the function
$$t\mapsto\int_0^t\!\!\int (e^{s\Delta}a_1e^{s\Delta}a_2)(x,s)\,dx\,ds.$$
This last problem is finally reduced to a linear system that can be solved
with elementary linear algebra.

\medskip
The spatial decay at infinity of the velocity
field is known to be closely related to special algebraic relations in terms of the moments
$\int x^\alpha \hbox{curl}(u)(x,t)\,dx$ of the vorticity~$\hbox{curl}(u)$
of the flow, see~\cite{Gallay02}.
Thus, one could restate the theorem in an equivalent way
in terms of identities between such moments
for different values of~$\alpha\in\N^d$, which are satisfied at the time $t^*_i$
but brake down when $t=t_i'$.

\subsection{A concentration effect of a different nature}

The concentration-diffusion effects described
in Theorem~\ref{theorem1}
genuinely depend on the very
special structure of the nonlinearity $\P\nabla\cdot (u\otimes u)$,
more than to the presence of the $\Delta u$ term.
Even though this result is not known for inviscid flows yet,
it can be expected that a similar property should be observed also
for the Euler equation.

On the other hand,
the Laplace operator, commonly associated with diffusion effects,
can be responsible also of concentration phenomena, of a different nature.
For example, it can happen that $a(x)$ is a \emph{non decaying\/}
(or very slowly decaying)
vector field, but such that the unique strong solution $u(x,t)$ of the Navier--Stokes
system
have a quite fast pointwise decay as $|x|\to\infty$ (say, $\sim|x|^{-d-1}$).
This is typically the case when $a$ has rapidly increasing oscillations in the far field.
We will  discuss this issue in Section~\ref{secKato}.
Though elementary, the examples of flows presented in that section have some interest,
being closely related to a problem addressed by Kato
 about strong solutions
in $L^p(\R^d)$ when $p<d$,
in his well known paper~\cite{Ka84}.

\subsection{Notations}
\label{secnot}

Troughout the paper, if $u=(u_1,\ldots,u_d)$ is a vector field with components in a
linear space $X$, we will write $u\in X$, instead of $u\in X^d$.
We will adopt a similar  convention for the tensors of the form $u\otimes u$.
We denote with $e^{t\Delta}$ the heat semigroup.

Let $B(0,1)$ be the unit ball in $\R^d$
and $\phi\in\mathcal{S}(\R^d)$ a function
satisfying
\begin{equation}
\label{condphi}
\begin{aligned}
&\widehat \phi\in C^\infty_0(\R^d),  \quad \hbox{supp}\,\widehat\phi\subset B(0,1),\quad \hbox{$\widehat \phi$ radial},
\quad \widehat\phi\ge0,
\quad\int|\widehat\phi|^2=1/d.
\end{aligned}
\end{equation}
Our definition for the Fourier transform is $\widehat\phi(\xi)=\int \phi(x)e^{-i\xi\cdot x}\,dx$.
Then we set
\begin{equation}
 \widehat \phi^\delta(\xi)=\frac{\widehat\phi(\xi/\delta)}{\delta^{d/2}},
\qquad\delta>0.
\end{equation}
Next we define the orthogonal transformation $\widetilde\;\,\colon\R^d\to\R^d$, by
\begin{equation}
 \label{orth}
\begin{aligned}
&\widetilde\alpha =(\alpha_2,\alpha_1), \qquad    &\hbox{if $\alpha=(\alpha_1,\alpha_2)\in\R^2$},\\
&\widetilde\alpha =(\alpha_2,\alpha_3,\alpha_1), \quad&\hbox{if $\alpha=(\alpha_1,\alpha_2,\alpha_3)\in\R^3$}.
\end{aligned}
\end{equation}

We define the $\hbox{curl}(\cdot)$ operator by
\begin{equation*}
\hbox{curl}\,\psi=(-\partial_2,\partial_1)\psi, \qquad\hbox{if $\psi\colon\R^2\to\R$}
\end{equation*}
and by
\begin{equation*}
 \hbox{curl}\,{\psi}=
\left(
\begin{matrix}
 \partial_2\psi_3-\partial_3\psi_2\\
\partial_3\psi_1-\partial_1\psi_3\\
\partial_1\psi_2-\partial_2\psi_1
\end{matrix}
\right)
\qquad
\hbox{if $\psi\colon\R^3\to\R^3$}.
\end{equation*}
The notation $f(x,t)=\mathcal{O}_t(|x|^{-\alpha})$ as $|x|\to\infty$ means that $f$
satisfies, for large~$|x|$, a bound of the form $|f(x,t)|\le A(t)|x|^{-\alpha}$, for some  function
$A$ locally bounded in $\R^+$.

We shall make use of the usual Kronecker symbol, $\delta_{j,k}=1$ or~$0$, if $j=k$ or $j\not=k$.


\section{Nonlinear concentration-diffusion effects}
\label{section2}

\subsection{The analyticity of the flow map}
\label{sec2.1}
In this subsection we recall a few well known facts.

Let  $B$ be the Navier--Stokes  bilinear operator, defined by
\begin{equation*}
 \label{Buv}
B(u,v)(t)\equiv-\int_0^t e^{(t-s)\Delta}\P\nabla\cdot(u\otimes v)(s)\,ds.
\end{equation*}
Then the Navier--Stokes equations can be written in the following integral form
\begin{equation}
 \label{NSI}
u=u_0+B(u,u), \qquad u_0=e^{t\Delta}a,\qquad \hbox{div}(a)=0.
\end{equation}

\medskip
Even though in the sequel we will only deal with ``concrete'' functional spaces,
the problematic is better understood in an abstract setting:
we will present it as formulated in the paper by P.~Auscher and Ph.~Tchamitchian
\cite{AuscT}.
Let $\mathcal{F}$ be a Banach space, $u_0\in\mathcal{F}$ and let
$B\colon\mathcal{F}\times\mathcal{F}\to\mathcal{F}$ be a continuous bilinear operator,
with operator norm~$\|B\|$.

Let us introduce the nonlinear operators $T_k\colon\mathcal{F}\to\mathcal{F}$, $k=1,2\ldots$, defined
through the formulae
\begin{equation*}
 \label{indu}
\begin{split}
 &T_1={\rm Id}_{\mathcal{F}}\\
& T_k(v)\equiv\sum_{l=1}^{k-1}B(T_l(v),T_{k-l}(v)), \quad k\ge2.
\end{split}
\end{equation*}
Then the following estimate holds:
\begin{equation*}
 \label{ATE}
\|T_k(u_0)\|_{\mathcal{F}}\le \frac{C}{\|B\|}k^{-3/2}\bigl(4\,\|B\|\,\,\|u_0\|_{\mathcal{F}}\bigr)^k.
\end{equation*}
Moreover, $T_k$
is the restriction to the diagonal of $\mathcal{F}^k= \mathcal{F}\times\cdots\times \mathcal{F}$
of a $k$-multilinear operator $\colon \mathcal{F}^k\to\mathcal{F}$.

Under the smallness assumption
\begin{equation}
 \label{absmall}
\|u_0\|_{\mathcal{F}}\le 1/(4\|B\|),
\end{equation}
the series
\begin{equation}
\label{snc}
\Psi(u_0)\equiv\sum_{k=1}^\infty T_k( u_0 ),
\end{equation}
is absolutely convergent in~$\mathcal{F}$
and its sum $\Psi(u_0)$ is a solution
of the equation
$ u=u_0+B(u,u)$.
Furthermore, $\Psi(u_0)$ is the only solution in the closed ball
$\overline{B_{\mathcal{F}}}(0,\frac{1}{2\|B\|})$.
The proof of these claims would be a straightforward application
of the contraction mapping theorem under the more restrictive condition
$\|u_0\|_{\mathcal{F}}<1/(4\|B\|)$. The proof of the more subtle version stated here can
be found in~\cite{AuscT}, \cite{Lem02}.

\medskip
Coming back to Navier--Stokes,
in the proof of Theorem~\ref{theorem1} we will need
to write the solution of the Navier--Stokes system as $u=\Phi(a)$, where
\begin{equation}
 \label{sols}
\Phi(a)(t)\equiv\sum_{k=1}^\infty T_k(e^{t\Delta}a),
\end{equation}
the series being absolutely convergent in $\mathcal{C}([0,\infty),L^2(\R^d))$.
There are several ways to achieve this, and  one of the simplest 
(which goes through in all dimension~$d\ge2$) is the following:
we choose $\mathcal{F}$ as the space of all functions~$f$ in $\mathcal{C}([0,\infty),L^2(\R^d))$,
such that $\|f\|_{\mathcal{F}}<\infty$, where
\begin{equation}
\label{newf}
\|f\|_{\mathcal{F}}\equiv  \hbox{ess\,sup}_{x,t}(1+|x|)^{d+1}|f(x,t)|\,\,\,+\,\,\,
\hbox{ess\,sup}_{x,t}(1+t)^{(d+1)/2}|f(x,t)|.
\end{equation}
The bicontinuity of the bilinear operator~$B$ is easily proved in
this space~$\mathcal{F}$ (see \cite{Miy02}).
Indeed, one can prove this only using the well known scaling relations and 
pointwise estimates on the kernel~$F(x,t)$
of the operator $e^{t\Delta}\P\hbox{div}$:
\begin{equation*}
F(x,t)=t^{-(d+1)/2}F(x/\sqrt t,1), \qquad |F(x,1)|\le C(1+|x|)^{-d-1}.
\end{equation*}
We can conclude that there is a constant $\eta_d>0$, depending only on~$d$,
such that if
\begin{equation}
\label{sett}
\|e^{t\Delta}a\|_{\mathcal{F}}<\eta_d
\end{equation}
then there is a solution $u=\Phi(a)\in \mathcal{F}$ of the Navier--Stokes equations
such that the series~\eqref{sols} is absolutely convergent in the~$\mathcal{F}$-norm.
The absolute convergence of such series in~$\mathcal{C}([0,\infty),L^2(\R^d))$
is then straightforward under the smallness assumption~\eqref{sett}.

The finiteness of~$\|e^{t\Delta}a\|_{\mathcal{F}}$ can be ensured, {\it e.g.\/}
by the two conditions $\int |a(x)|(1+|x|)\,dx<\infty$ and 
$\hbox{ess\,sup}_{x\in\R^d} (1+|x|)^{d+1}|a(x)|<\infty$.
The smallness condition~\eqref{sett} could be be slightly relaxed, see~\cite{BraV07}.

\subsection{The construction of the initial datum}

This section devoted to a constructive proof of the following Lemma.

\begin{lemma}
\label{prop1}
Let $d=2,3$,  and $\epsilon>0$. Let also $N\in\N$ and $0<t_1<\cdots<t_N$ be a finite sequence.
Then there exists a divergence-free vector field $a=(a_1,\ldots,a_d)\in\mathcal{S}(\R^d)$,
such that
\begin{equation}
 \label{inav}
\widetilde a(x)=a(\widetilde x), \qquad x\in\R^d
\end{equation}
(see Section~\ref{secnot} for the notations)
and such that the function
$E(a)(t)\colon\R^+\to\R$, defined by
\begin{equation}
 \label{Eap}
E(a)(t)\equiv  - \int_0^t\!\!\int e^{s\Delta}a_1(x)\,e^{s\Delta}a_2(x)\,dx\,ds,
\end{equation}
changes sign inside $(t_i-\epsilon,t_i+\epsilon)$, for $i=1,\ldots,N$.
\end{lemma}

\Proof
It is convenient to separate the two and three-dimensional cases

\medskip
\noindent
{\bf The case $\boldsymbol{d=2}$.}
We start setting, for each~$\alpha\in\R^2$,
\begin{equation}
 \label{psia}
\widehat \psi_\alpha(\xi)\equiv \widehat\phi(\xi-\alpha)+\widehat \phi(\xi+\alpha)-\widehat\phi(\xi-\widetilde\alpha)
	-\widehat\phi(\xi+\widetilde\alpha),
\end{equation}
for some $\phi\in\mathcal{S}(\R^2)$ satisfying conditions~\eqref{condphi}.
Next we introduce the divergence-free vector field $a_\alpha(x)$, through the relation
\begin{equation}
 \label{aalpha}
\widehat a_\alpha(\xi)=\biggl(
\begin{aligned}
 -i\xi_2\widehat \psi_\alpha(\xi)\\ 
   i\xi_1\widehat \psi_\alpha(\xi)
\end{aligned}
\biggr).
\end{equation}
Note that $a_\alpha\in\mathcal{S}(\R^2)$ is real-valued
(because $\widehat\psi_\alpha$ is real-valued and such that $\widehat\psi_\alpha(\xi)=\widehat\psi_\alpha(-\xi)$)
and satisfies the fundamental symmetry condition
\begin{equation}
 \label{syma}
\widetilde a_\alpha(x)=a_\alpha(\widetilde x).
\end{equation}
Next we  define $\widehat\psi^\delta$, $a_\alpha^\delta$ as before,
by simply replacing $\widehat\phi$ with $\widehat\phi^\delta$ in the corresponding definitions.
The vector fields $a^\delta_\alpha$ will be our ``building blocks'' of our initial datum.

Applying the Plancherel theorem to the right-hand side of~\eqref{Eap}
and using the symmetry relations $\widehat\psi_\alpha(\xi)=\widehat\psi_\alpha(-\xi)$,
$|\widehat\psi_\alpha(\xi)|=|\widehat\psi_\alpha(\tilde\xi)|$,
we get
\begin{equation}
\label{Eaa}
\begin{split}
E(a_\alpha)(t)
&=\int (1-e^{-2t|\xi|^2}) \frac{\xi_1\xi_2}{2|\xi|^2} |\widehat\psi_\alpha(\xi)|^2\,d\xi\\
&=2\int_{\xi_1\ge|\xi_2|} (1-e^{-2t|\xi|^2}) \frac{\xi_1\xi_2}{|\xi|^2} |\widehat\psi_\alpha(\xi)|^2\,d\xi.
\end{split}
\end{equation}

\medskip
From now on, the components of $\alpha\in\R^2$
will be assumed to satisfy the following conditions
\begin{equation}
 \label{assal}
\begin{cases}
 \alpha_1>|\alpha_2|\\ \alpha_2\not=0.
\end{cases}
\end{equation}
This guarantees that for a sufficiently small $\delta>0$
({\it i.e.\/} when $\alpha_1>|\alpha_2|+\delta\sqrt 2$),
we have
\begin{equation}
\label{twen}
E(a_\alpha^\delta)(t)=2\int (1-e^{-2t|\xi|^2})\frac{\xi_1\xi_2}{|\xi|^2} |\widehat \phi^\delta(\xi-\alpha)|^2\,d\xi.
\end{equation}

\medskip
If we set
\begin{equation}
\label{Eapp} 
E^{\rm app}_\alpha(t)\equiv \bigl(1-e^{-2t|\alpha|^2}\bigr)\frac{\alpha_1\alpha_2}{|\alpha|^2}
\end{equation}
then we immediately obtain
\begin{equation*}
 E(a^\delta_\alpha)(t)\to E^{\rm app}_\alpha(t), \qquad\hbox{as $\delta\to0$}
\end{equation*}
uniformly with respect to~$t\ge0$.

\medskip

We now associate to  $(t_1,\ldots,t_N)$
two  more sequences (to be chosen later)
$(\lambda_1,\ldots,\lambda_{N+1})\in\R_+^{N+1}$ and 
$(\alpha_1,\ldots,\alpha_{N+1})\in\R^{2(N+1)}$,
where $\alpha_j=(\alpha_{j,1},\alpha_{j,2})$.
First we require that the components $\alpha_{j,1}$ and $\alpha_{j,2}$
satisfy  condition~\eqref{assal} for all $j=1,\ldots,N+1$ and
that $\alpha_j\not=\alpha_{j'}$, for $j\not=j'$ and $j,j'=1,\ldots,N+1$.
This  second requirement
ensures that the supports of $\widehat a^\delta_{\alpha_j}$ and
$\widehat a^\delta_{\alpha_{j'}}$ are disjoint when 
$\delta$ is sufficiently small.

\medskip
We now consider the initial data of the form
\begin{equation}
 \label{init}
a^\delta(x)\equiv\sum_{j=1}^{N+1}\lambda_j a^\delta_{ \alpha_j}(x).
\end{equation}
Owing to the condition on the supports of $\widehat a^\delta_{\alpha_j}$, 
we see that for $\delta>0$ small enough,
\begin{equation}
 \label{Ead}
E(a^\delta)(t)=\sum_{j=1}^{N+1} \lambda_j^2 E(a^\delta_{\alpha_j})(t).
\end{equation}
Let
\begin{equation}
\label{lamm} 
  \mu_j=\frac{\lambda_j^2\,\alpha_{j,1}\,\alpha_{j,2}}{ |{\alpha_j}|^2 }
\qquad\hbox{and}\qquad
A_j=e^{-2|{\alpha}_j|^2}.
\end{equation}
Thus, as $\delta\to0$, we get
$E(a^\delta)(t)\to E^{\rm app}(t)$,
uniformly in $[0,\infty)$,
where
\begin{equation}
 E^{\rm app}(t)=\sum_{j=1}^{N+1} \mu_j\bigl(1-A_j^t\bigr).
\end{equation}
Let us observe that
\begin{equation}
\label{deE} 
\frac{{\rm d} E^{\rm app}}{{\rm d} t}(t)=
 -\sum_{j=1}^{N+1}\mu_j \,\log(A_j) A_j^t.
\end{equation}
We want to determine $(\lambda_{1},\ldots,\lambda_{N+1})$ and $(\alpha_1,\ldots,\alpha_{N+1})$
in a such way that $E^{\rm app}(t)$ vanishes at $t_1,\ldots,t_N$,
and changing  sign at those points.
This leads us to study the system of~$N$ equalities and~$N$ `non-equalities',
\begin{equation}
 \label{deri}
\left\{
\begin{aligned}
 &E^{\rm app}(t_i)=0\\
&\frac{{\rm d}E^{\rm app}}{{\rm d}t}(t_i)\not=0.\\
\end{aligned}
\right.
\qquad i=1,\ldots,N.
\end{equation}

Let us choose $|\alpha_j|^2=\gamma j$, for an arbitrary $\gamma>0$
 (this choice is not essential, but will greatly simplify the calculations)
and set $T_i=e^{-2\gamma t_i}$.
Recalling~\eqref{lamm}, we get $A_j^{t_i}=T_i^{j}$ and $\log(A_j)=-2\gamma j$.
In order to study the system~\eqref{deri}, we introduce the $(N+1)^2$-matrix
\begin{equation}
\label{matr}
M:=
\left(
 \begin{matrix}
  1-T_1& 1-T_1^2& \cdots &1-T_1^{N+1} \\
\vdots  & \vdots &              & \vdots \\
1-T_N & 1-T_N^2 &\cdots &1-T_{N}^{N+1}\phantom{\Bigl|} \\
T_1 & 2T_1^2      &\cdots & (N+1)T_1^{N+1}
 \end{matrix}
\right)
\end{equation}
We claim that $\det{M}\not=0$.
Indeed, by an explicit computation,
\begin{equation*}
\det(M)= - T_1\,(1-T_1)\prod_{i=1}^N(1-T_i)\prod_{i=2}^N (T_1-T_i)\prod_{1\le i<i'\le N} (T_{i'}-T_i).
\end{equation*}
Recalling that $T_i=e^{-2\gamma t_i}\in(0,1)$ and that $t_i\not=t_{i'}$ proves our claim.
The above formula can be checked by induction.
Otherwise, one can reduce~$M$ after elementary factorizations
to a Vandermonde-type matrix
(see \cite{Kratt1} for explicit formulae 
on determinants).

Then, for any $c\not=0$,  the linear system with unknown $\boldsymbol{\mu}=(\mu_1,\ldots,\mu_{N+1})$,
\begin{equation}
 \label{obl}
M \boldsymbol{\mu}=
\left(
\begin{matrix}
0\\ \vdots\\ 0\\ c 
\end{matrix}
\right)
\end{equation}
has a unique solution~$\boldsymbol{\mu^*}\in\R^{N+1}$, $\boldsymbol{\mu^*}\not=0$.
By our construction, the function~$E^{\rm app}(t)$ obtained taking $\boldsymbol{\mu}=\boldsymbol{\mu^*}$
satisfies the $N$-equations and the first `non-equality' of the system~\eqref{deri}.
More precisely, we get $({\rm d}/{\rm d}t)E^{\rm app}(t_1)=c/(2\gamma)\not=0$.
The other $N-1$ `non-equalities' of the system~\eqref{deri} are then automatically fulfilled.
Indeed if, otherwise, we had $({\rm d}/{\rm d}t)E^{\rm app}(t_i)=0$,
for some $i=2,\ldots,N$, then the matrix obtained replacing in~\eqref{matr}
the last line with 
$$\begin{matrix}T_i \,&\, 2T_i^2  \,     &\,\cdots \,&\, (N+1)T_i^{N+1}\end{matrix}$$
would have been of determinant zero, thus contredicting our preceding formula for $\det(M)$.

\medskip
By conditions~\eqref{lamm}, for all $j=1,\ldots,N+1$, the real number $\mu^*_j$ defines
(in a non-unique way)  a real $\lambda_{j}^*$ and a vector $(\alpha^*_{j,1},\alpha^*_{j,2})$
with components satisfying~\eqref{assal}, such that $(\lambda_{1},\ldots,\lambda_{N+1})\not=(0,\ldots,0)$.

We now consider the initial data $a^\delta_*$ obtained from formula~\eqref{init},
choosing $\lambda_j=\lambda_j^*$ and $\alpha_j=\alpha_j^*$
for  $j=1,\ldots,N+1$.
With this choice, the corresponding function
$E^{\rm app}(t)$ changes sign in each neighborhood of~$t_i$.

By the uniform convergence of $E(a^\delta_*)(t)$ to $E^{\rm app}(t)$ as $\delta\to0$,
we see that if $\delta>0$ is small enough then
$E(a^\delta_*)$ changes sign in the interval $(t_i-\epsilon,t_i+\epsilon)$, 
for $i=1,\ldots,N$.

The conclusion of Lemma~\ref{prop1} in the two dimensional case now follows.

\medskip\noindent
{\bf The case $\boldsymbol{d=3}$.}
We only indicate the modifications that have to be done to the above arguments.
Let us set, for $\alpha\in\R^3$,
\begin{equation}
 \label{3dd}
\widehat a_{\alpha}(\xi)=
\left(
\begin{matrix} (i\xi_2-i\xi_3)\widehat\psi_\alpha(\xi)\\
 (i\xi_3-i\xi_1)\widehat\psi_\alpha(\xi)\\
 (i\xi_1-i\xi_2)\widehat\psi_\alpha(\xi)
\end{matrix}
\right)
\end{equation}
where $\psi$
is defined by
\begin{equation*}
 \begin{aligned}
  \widehat\psi_\alpha(\xi)= &\widehat\phi(\xi-\alpha)+\widehat\phi(\xi-\widetilde\alpha)
+\widehat\phi(\xi-\widetilde{\widetilde\alpha})\\
&\qquad
 +\widehat\phi(\xi+\alpha)+\widehat\phi(\xi+\widetilde\alpha)+\widehat\phi(\xi+\widetilde{\widetilde{\alpha}}).
 \end{aligned}
\end{equation*}
In this way, $a\in\mathcal{S}(\R^3)$ is a real valued divergence-free vector field
with the rotational symmetry~\eqref{inav}. As before, we can define also the rescaled
vector field $a^\delta$, for any $\delta>0$.

The conditions to be imposed on the components of~$\alpha=(\alpha_1,\alpha_2,\alpha_3)$ are now
\begin{equation}
 \label{assal3}
\min(\alpha_2,\alpha_3)>\alpha_1^+, \qquad \alpha_2\not=\alpha_3,
\end{equation}
where $\alpha_1^+=\max(\alpha_1,0)$.
Geometrically, the inequality in~\eqref{assal3} corresponds to cutting~$\R^3$ into
six congruent regions, that can be obtained from each other through the orthogonal
transforms $\alpha\mapsto\widetilde\alpha$ and $\alpha\mapsto -\alpha$,
and then selecting one of these regions.
If $\delta>0$ is small enough, in a such way that
$\alpha+B(0,\delta)$ is contained in the region~$\Gamma\subset\R^3$ defined by~\eqref{assal3},
then we get, recalling the definition~\eqref{Eap} of~$E(a_\alpha)$, the three-dimensional counterpart of~\eqref{twen}:
\begin{equation*}
 E(a_\alpha)(t)
=3\int(1-e^{-2t|\xi|^2})\frac{(\xi_1-\xi_3)(\xi_2-\xi_3)}{|\xi|^2}|\widehat\phi^\delta(\xi-\alpha)|^2\,d\xi.
\end{equation*}
(Here one first applies Plancherel Theorem, then the integral over $\xi\in\R^3$
is rewritten, because of  the symmetries of $\widehat a_\alpha$, as $6\int_{\Gamma}\dots\,d\xi$).
We have, for all $t\ge0$ and as $\delta\to0$,
\begin{equation*}
E(a^\delta_\alpha)(t)\to E^{\rm app}_\alpha(t)
\equiv (1-e^{-2t|\alpha|^2})\frac{(\alpha_1-\alpha_3)(\alpha_2-\alpha_3)}{|\alpha|^2}.
\end{equation*}
For $\alpha_j=(\alpha_{j,1},\alpha_{j,2},\alpha_{j,3})\in\R^3$, $j=1,\ldots,N+1$, with components satisfying~\eqref{assal3},
we now set
\begin{equation}
\label{lamm3} 
  \mu_j=\frac{\lambda_j^2\,(\alpha_{j,1}-\alpha_{j,3})(\alpha_{j,2}-\alpha_{j,3})}{ |{\alpha_j}|^2 }.
\end{equation}
As in the two-dimensional case, it is possible to choose the phases~$\alpha_j$ in a such way that
$|\alpha_j|^2=\gamma j$, where $\gamma>0$ is arbitrary. This leaves unchanged the definitions of $A_j$ and~$T_i$.

Now solving, as before, the linear system~\eqref{obl}, shows that
it is possible
to construct an initial datum
of the form
$$ a^\delta=\sum_{j=1}^{N+1} \lambda_j a_{\alpha_j}^\delta,$$
such that $E(a)(t)$ has a non-constant sign inside the intervals $(t_i-\epsilon,t_i+\epsilon)$,
for
$i=1,\ldots,N$.

Lemma~\ref{prop1} is now established.

\nobreak
\endProof

\subsection{End of the proof of Theorem~\ref{theorem1}}

We are now in position to deduce from Lemma~\ref{prop1}
and the facts recalled in Section~\ref{sec2.1} the conclusion of Theorem~\ref{theorem1}.

\medskip
\noindent
{\bf Step 1.}
\emph{Constructing a solution such that $\,t\mapsto \int_0^t\!\!\int u_1u_2(x,s)\,dx\,ds$
changes sign near~$t_1$, $t_2$,\dots,$t_N$.}

\medskip
\noindent
Let us consider the initial datum $a$ constructed in Lemma~\ref{prop1}.
If necessary, we modify~$a$ by multipliying it by a small constant
$\eta_0>0$ in order to ensure that the corresponding solution~$u$ of the Navier--Stokes system
is defined globally in time.
Without loss of generality, we can and do assume $\eta_0=1$.

With the same notation of Section~\ref{sec2.1}, let
$K(a)\colon\R^+\to\R^{d\times d}$,
\begin{equation}
K(a)(t)\equiv\int_0^t\!\!\int \bigl(\Phi(a)\otimes\Phi(a)\bigr)(x,s)\,dx\,ds.
\end{equation}
Let~$\mathcal{F}$ be the space defined in Section~\ref{sec2.1}, and normed by~\eqref{newf}.
For $\eta>0$ small enough we can write, by expansion~\eqref{sols},
\begin{equation}
 \label{snd}
K(\eta\, a)(t)=\sum_{k=2}^\infty \eta^k\,S_{k}(e^{t\Delta}a)(t)
\end{equation}
where $S_{k}\colon \mathcal{F}\to\R^{d\times d}$
is the restriction to the diagonal of $\mathcal{F}^k=\mathcal{F}\times\cdots\times \mathcal{F}$
of a $k$-multilinear operator defined  from
$\mathcal{F}^k$ to $\R^{d\times d}$.

The embedding of~$\mathcal{F}\subset \mathcal{C}([0,\infty),L^2(\R^d))$
implies that, for all~$t>0$, the series in~\eqref{snd} is absolutely convergent.
But,
\begin{equation*}
 \label{S2}
S_2(e^{t\Delta}a)(t)=\int_0^t\!\!\int \bigr(e^{s\Delta}a\otimes e^{s\Delta}a\bigr)(x,s)\,dx\,ds,
\end{equation*}
so that, according to the notations of Lemma~\ref{prop1},
\begin{equation*}
 \label{SE}
\bigl(S_2(a)\bigr)_{1,2}(t)=-E(a)(t).
\end{equation*}
For any time $t$ such that $E(a)(t)>0$ (respectively, $E(a)(t)<0$), we can find a small 
$\eta_{t}>0$ such that the solution $\Phi(\eta\, a)$ of the Navier--Stokes system
starting from $\eta \, a$ satisfies, for all $0<\eta<\eta_t$, 
\begin{equation*}
K(\eta\, a)_{1,2}(t)<0, \qquad\hbox{(respectively, $K(\eta\,a)_{1,2}(t)>0$)}.
\end{equation*}
By Lemma~\ref{prop1},
this observation can be applied
for~$t$ belonging to two suitable finite sequences of times
(that can be taken arbitrarily close to $(t_1,\ldots,t_N)$ in the $\R^N$-norm).
We conclude that if $\eta>0$ is small enough then $K(\eta\,a)_{1,2}(t)$ changes sign
in $(t_i-\epsilon,t_i+\epsilon)$, for $i=1,\ldots,N$.
In particular, because of the continuity of the map $t\mapsto K(\eta\, a)(t)$,
there is a point $t^*_i$ inside each one of these intervals where $K(\eta\,a)_{1,2}(t^*_i)=0$.

\medskip
\noindent
{\bf Step 2.}
\emph{Analysis of the orthogonality relations $\int_0^t\!\!\int (u_ju_k)(x,s)\,dx\,ds=c(t)\delta_{j,k}$}.

\medskip
\noindent
Let $u=\Phi(\eta\, a)$ be the solution constructed in the first step.
Notice that
$\widetilde u(x,t)=u(\widetilde x,t)$ for all $\hbox{$t\ge0$}$,
{\it i.e.\/}, the condition $\widetilde a(x)=a(\widetilde x)$ propagates during the evolution.
Indeed, this is simple consequence of the invariance of the Navier--Stokes equations
under the tranformations of the orthogonal group $O(d)$ and the uniqueness of strong solutions.
Thus, $\int u_1^2\,dx=\int u_2^2\,dx$. 
When $d=3$, such integrals equal, of course, $\int u_3^2\,dx$, and we have also $\int u_1u_2\,dx=\int u_2u_3\,dx=\int u_3u_1\,dx$.
We deduce that, for all~$t\ge0$ (we denote here by $I$ the $d\times d$ identity matrix):
\begin{equation}
\label{iff1}
 \exists \, c(t)\colon\; K(\eta\,a)(t)=c(t)I
\qquad\hbox{if and only if}\qquad
K(\eta\,a)_{1,2}(t)=0.
\end{equation}
(This equivalence is no longer valid for $d\ge4$. For a proof of the theorem in
the higher dimensional case one should consider flows invariant under
larger discrete subgroups of $O(d)$).

\medskip
\noindent
{\bf Step 3.}
\emph{The far-field asymptotics of the velocity field.}

\medskip\noindent
By the result of \cite{BraV07} (see Theorem~1.2 or Theorem~1.7,
applied in the particular case of the datum $\eta\,a\in\mathcal{S}(\R^d)$)
we know that, for all $t>0$,
\begin{equation}
\label{bva}
 u(x,t)=\nabla_x\Pi(x,t)+\mathcal{O}_t(|x|^{-d-2}), \qquad \hbox{as $|x|\to\infty$},
\end{equation}
where,
\begin{equation}
\label{qua}
 \Pi(x,t)=\gamma_d\sum_{h,k}\biggl(\frac{\delta_{h,k}}{d|x|^d}-\frac{x_hx_k}{|x|^{d+2}}\biggr)\cdot K(\eta\,a)_{h,k}(t)
\end{equation}
and $\gamma_d\not=0$ is a constant.

Concerning the first term $\nabla_x\Pi$ on the right-hand side of~\eqref{bva},
for each fixed~$t>0$ two situations can occur.
Either the function $x\mapsto \nabla_x\Pi(x,t)$ is identically zero,
or this function is homogeneous of degree exactly~$-d-1$.
But, for all fixed $t>0$ (see \cite{BraV07}, Proposition~1.6),
\begin{equation}
 \label{iff2}
\nabla_x\Pi(\cdot,t)\equiv0\qquad\hbox{if and only if} \qquad \exists \,c(t)\in\R\colon\; K(\eta\,a)_{h,k}(t)=c(t)\delta_{h,k}.
\end{equation}

Combining conditions~\eqref{iff1}-\eqref{iff2} with the asymptotic profile~\eqref{bva},
we deduce from the analysis we made in Step~1,
the upper bound
$$|u(x,t^*_i)|\le C|x|^{-d-2}, \qquad\hbox{$i=1,\ldots,N$} $$
and the lower bound,
$$|u_j(x,t'_i)|\ge c_{\omega} |x|^{-d-1}, \qquad\hbox{$i=1,\ldots,N$}, \quad j=1,\ldots,d,$$
for all $|x|$ large enough 
and some points $t'_i$ distant less than~$\epsilon$ from $t_i$.
Here, $C>0$ is independent on~$x$ and
$c_\omega$ is independent on~$|x|$, but will depend on the projection
$\omega=\frac{x}{|x|}$ of~$x$ on the sphere~$\S^{d-1}$.
In fact, we can take $c_{\omega}>0$, unless $\partial_{x_j}\Pi(\omega,t'_i)$ has a zero
at the point~$\omega\in\S^{d-1}$.
But one deduces from~\eqref{qua} that the zeros of $\partial_{x_j}\Pi(\omega,t'_i)$
are exactly the zeros on the unit sphere of a homogeneous polynomial
of degree three. Therefore, $c_\omega>0$ for almost every $\omega\in\S^{d-1}$.

Theorem~\ref{theorem1} is now established.

\endProof

\subsection{Explicit examples}
Let us exhibit an explicit example of our construction,
in the simplest case $N=1$.
We set, for $d=2$, and a function~$\phi$ satisfying conditions~\eqref{condphi},
\begin{equation*}
\begin{split} 
a(x)=
\eta\left(
\begin{matrix}
 -\partial_2\\ \partial_1
\end{matrix}
\right)
\biggl[\phi^\delta(x)
\Bigl( &\sqrt 3\cos(\sqrt 3x_1+x_2)-\sqrt 3\cos(x_1+\sqrt 3 x_2) \\
&+\sqrt 2\cos(\sqrt 6x_1-\sqrt 2x_2)-\sqrt 2\cos(-\sqrt 2x_1+\sqrt 6x_2)\Bigr)\biggr].
\end{split}
\end{equation*}
This expression is obtained taking, in Lemma~\ref{prop1},
$T_1=\frac{1}{2}$,  $\alpha_1=(\sqrt 3,1)$, $\alpha_2=(\sqrt 6,-\sqrt 2)$
and $\gamma=4$, and observing that
$\lambda_1^2=\frac{3}{2}\lambda_2^2$.
If $|\eta|$ and $\delta$ are positive and small enough, then the solution $u(x,t)$ of Navier--Stokes starting
from~$a$ concentrates/diffuses  when $t\simeq \frac{1}{8}\log(2)$.

\medskip
Three-dimensional examples are obtained in a very similar way:
choose again $T_1=\frac{1}{2}$, next $\alpha_1=(0,1,\sqrt 3)$, $\alpha_2=(0,\sqrt 6,\sqrt 2)$.
The relation between the coefficients $\lambda_1$ and $\lambda_2$
is now $\lambda_1^2=\frac{\sqrt 3}{2}\lambda_2^2$.
This leads us to introduce the function
\begin{subequations}
\begin{equation}
\label{sub1}
\begin{split}
f^\delta(x)\equiv
\phi^\delta(x)\biggl[
 &3^{1/4}\Bigl(\cos(x_2+\sqrt 3 x_3)+\cos(x_1+\sqrt 3 x_2)+\cos(\sqrt 3 x_1+x_3)\Bigr)\\
&
+\sqrt 2\Bigl(\cos(\sqrt 6 x_2+\sqrt 2 x_3)+\cos(\sqrt 6 x_1+\sqrt 2 x_2)+\cos(\sqrt 2 x_1+\sqrt 6 x_3)\Bigr)\biggr],
\end{split}
\end{equation}
where~$\delta>0$ and  $\phi$ satisfies conditions~\eqref{condphi}.
If~$a$ is the vector field
\begin{equation}
\label{sub2} a=\eta\, {\rm curl}(f^\delta,f^\delta,f^\delta),
\end{equation}
\end{subequations}
with $|\eta|,\delta>0$ small enough,
then the solution arising from~$a$ concentrates/diffuses, as before, around~$t\simeq \frac{1}{8}\log(2)$.

\begin{remark}
The smallness of~$\eta$ was important for the applicability of our analyticity argument.
However, it would be interesting to know whether three-dimensional data as those constructed in~\eqref{sub1}-\eqref{sub2},
with~\emph{large} coefficients~$\eta$,
still feature some kind of concentration effects in finite time.
\end{remark}


\section{Linear concentration-diffusion effects}
\label{secKato}

In this section we would like to give an elementary example
of a Navier--Stokes flow featuring a concentration effect of a quite different
nature. A similar effect can be observed if we replace the Navier--Stokes nonlinearity
by a more general term (not necessarily quadratic).
Let us focus, for example, on the three-dimensional case.

Consider the two initial data in~$\R^3$ (here $\eta\not=0$ is a constant)
\begin{equation}
 \label{da1}
\dot a(x)=\eta\Bigl( -\partial_{2}\bigl[ \log(e+|x|^2)^{-1}\bigr]\,,\,
\partial_1\bigl[\log(e+|x|^2)^{-1}\bigr]\,,\, 0 \Bigr)
\end{equation}
and
\begin{equation}
 \label{da2}
\overline a(x)=\dot a(x)\sin(|x|^2).
\end{equation}
In both cases, a direct computation shows that such data are divergence-free and belong to $L^3(\R^3)$,
but not to $L^p(\R^3)$ for $1\le p<3$.

\medskip
In his well-known paper~\cite{Ka84}, Kato could obtain,
for a solution~$u$ of the Navier--Stokes system emanating from~$a$,
local and global results implying $u(t)\in L^q(\R^3)$ 
for $t>0$, with $q\ge3$, under the assumption $a\in L^3(\R^3)$.
He also obtained results implying $u(t)\in L^p\cap L^q(\R^3)$,
always for $q\ge3$,  under the stronger assumption $a\in L^p\cap L^3(\R^3)$,
with $1\le p<3$.
Immediately after stating his theorems, Kato observed the following

\begin{remark} {\bf (see \cite{Ka84}) }
``The spatial decay expressed by the property $u(t)\in L^q(\R^3)$
is of interest.
Note that $q$ is restricted to $q\ge3$.
\emph{We were able to give no results for $q<3$ under the {\rm [only]} assumption $a\in L^3(\R^3)$}''.
\end{remark}

The solutions arising from the data $\dot a$ and $\overline{a}$ (both uniquely defined in
$\mathcal{C}([0,T),L^3(\R^3))$, for some $0<T\le\infty$ {\it a priori\/} depending on $\eta$)
show that when one assumes only $a\in L^3(\R^3)$ then ``everything can happen''.
Indeed, it is not difficult to prove the following claims:

\vbox{
\begin{enumerate}
\item
The solution $u(t)$ arising from the datum~$\dot a$ \emph{does not belong to $L^q(\R^3)$\/},
for any $t\in(0,T)$ and any $q\in[1,3)$.
\item
The solution $u(t)$ arising from the datum~$\overline a$ \emph{does belong to $L^q(\R^3)$\/}
for all $t\in(0,T)$ and all $q\in[1,3)$.
\end{enumerate}
}

Therefore, in the latter case, the solution enjoys some kind of spatial concentration effect.
However this effect does not rely on special geometric features of the flow,
but only on the oscillatory character of~$\overline a$.

\medskip

The first claim is immediate, because whenever $|a(x)|\le C(1+|x|)^{-1}$ and $|\nabla a(x)|\le C(1+|x|)^{-2}$,
then arguing as in~\cite[Sec.~4]{Bra08}
one obtains for $t>0$ an estimate of the form $|e^{t\Delta}a(x)-a(x)|\le C(t)(1+|x|)^{-2}$.
Then, from  the equation $u=e^{t\Delta}a +B(u,u)$ (the notation is as in Section~\ref{sec2.1}),
the  estimate $|u(x,t)-a(x)|\le \bar C(t)(1+|x|)^{-2}$.
Applying this observation to~$\dot a$ yields the conclusion.

\medskip
Let us sketch the proof of our second claim.
The faster and faster oscillations of $\overline a$ imply that,
for $t>0$ and $|x|$ large enough,
$$ |e^{t\Delta}\overline a(x)| \le C_m(t)|x|^{-m-1}\log^{-2}|x|,
  \qquad m=0,1,2\ldots$$
where the functions $C_m(t)$ are locally bounded in $(0,\infty)$ and
satisfy $C_m(t)\sim t^{-m/2}$ as $t\to0$.
In particular, $\|e^{t\Delta}\overline a\|_{3/2}\le C_1(t)$.
Let us iterate the integral Navier--Stokes equation:
\begin{equation}
 \label{bii}
u(t)=e^{t\Delta}\overline a+B(e^{t\Delta} \overline a,e^{t\Delta} \overline a)
+2B(e^{t\Delta}\overline a,B(u,u))+B(B(u,u),B(u,u)).
\end{equation}
Now applying elementary H\"older and Young inequalities
in the right-hand side of~\eqref{bii} (we can freely use here that $\|u(t)\|_3$ is bounded in $(0,T)$),
we get, for all $t\in(0,T)$,
$\|u(t)\|_1\le C(t)$ with $C(t)\sim t^{-1}$ as $t\to0$.
By interpolation one obtains that $u(t)\in L^q(\R^3)$ in $(0,T)$ for $1\le p<3$, with
an estimate on the blow-up of the $L^q$-norm as $t\to0$.

Let us observe that for $0<t<T$,
even though $u(t)$ decays pointwise much faster than $\overline a$, as $|x|\to\infty$,
the decay of $u(t)$ does not exceed that of $|x|^{-4}$, accordingly with the
limitations on the spatial localization described in the introduction.

\bibliographystyle{amsplain}

\end{document}